\documentclass[12pt,leqno,twoside]{amsart}

\usepackage[latin1]{inputenc}
\usepackage[T1]{fontenc}
\usepackage[colorlinks=true, pdfstartview=FitV, linkcolor=blue, citecolor=blue, urlcolor=blue]{hyperref}
\usepackage{amstext,amsmath,amscd, bezier,indentfirst,amsthm,amsgen,enumerate, geometry}
\usepackage[all,knot,arc,import,poly]{xy}
\usepackage{amsfonts,color, soul}  
\usepackage{amssymb}
\usepackage{latexsym}
\usepackage{epsfig}
\usepackage{graphicx}
\usepackage{srcltx}
\usepackage{enumitem}

\usepackage{tikz-cd}
\usepackage[all]{xy}

\usepackage{tikz} 
\usetikzlibrary{through} 
\usetikzlibrary{patterns} 
\usetikzlibrary{intersections} 
\usetikzlibrary{matrix} 
\usetikzlibrary{cd} 

\newtheorem{theorem}{Theorem}[section]
\newtheorem{corollary}[theorem]{Corollary}
\newtheorem{lemma}[theorem]{Lemma}
\newtheorem{proposition}[theorem]{Proposition}
\theoremstyle{definition}
\newtheorem{definition}[theorem]{Definition}
\theoremstyle{remark}
\newtheorem{remark}[theorem]{\sc Remark}
\newtheorem{example}[theorem]{\sc Example}

\theoremstyle{claim}

\newcommand{\Sing}{{\rm{Sing\hspace{2pt}}}}

\newcommand{\rank}{{\rm{rank\hspace{2pt}}}}

\newcommand{\Disc}{{\rm{Disc\hspace{2pt}}}}

\newcommand{\e}{\varepsilon}
\newcommand{\dd}{{\rm{d}}}
\newcommand{\m}{\setminus}

\newcommand{\codim}{{\rm{codim\hspace{2pt}}}}

\newcommand{\bR}{{\mathbb R}}
\newcommand{\bC}{{\mathbb C}}

\begin{document}
	\title[Topology of first integrals via Milnor fibrations II]{Topology of first integrals via Milnor fibrations II}

	\author[Fernando Reis]{F. Reis}
	\address{Universidade Federal do Esp\'{i}rito Santo- UFES, Brazil}
	\email{fernando.reis@ufes.br}

	\author[Maico Ribeiro]{M. Ribeiro}
	\address{Universidade Federal do Esp\'{i}rito Santo - UFES, Brazil}
	\email{maico.ribeiro@ufes.br}

	\author[Eurípedes da Silva]{E. da Silva}
	\address{Instituto Federal do Cear\'{a} - IFCE, Brazil}
	\email{euripedes.carvalho@ifce.edu.br}
	
	\begin{abstract}
This survey is the continuation of a series of works aimed at applying tools from Singularity Theory to Differential Equations. More precisely, we utilize the powerfull Milnor's Fibration Theory to give geometric-topological classifications of first integrals of differential systems. In the previous paper, systems of first-order quasilinear partial differential equations were examined, focusing on the case of an isolated singularity. Now, we address both cases of isolated and \textit{non-isolated singularities} for more general dynamical systems (namely, \textit{foliations}) that admit at least one first integral. For this, we utilize recently established connections between harmonic morphisms and Milnor fibrations to provide topological and geometric descriptions of the foliations under consideration. In particular, we apply these results to analyze the graph of solutions of some quasilinear systems.

\vspace{0.5 cm}

Harmonic morphism, Milnor fibrations, Minimal submanifolds.


\vspace{0.5 cm}

	\end{abstract}
	
	\maketitle

	\section{Introduction}\label{section: intro} 
	
	In \cite{RRE2}, a discussion on the topological classification of solutions to systems of quasi-linear partial differential equations  was introduced, through the Singularity Theory. Beginning with a concise overview on the \textit{Integrability problem}, the work focuses on analyzing a system of first-order quasi-linear partial differential equations, given by
	\begin{equation}\label{eq: system 1}
		\sum_{i=1}^{n}a_i^{\lambda}(x,u)\frac{\partial u}{\partial x_{i}}= b^{\lambda}(x,u),
	\end{equation}
	where $u$ is one unknown function in $n$ independent variables $x=(x_1,\dots,x_n)$,  $\lambda =1,\dots,q$, with $1\leq q\leq n$ and  $a_i^{\lambda},b^{\lambda}$ are analytic functions defined on an open set $\Omega \subset \mathbb{R}^{n+1}$. Then, utilizing one of the fibered structures known as \textit{Milnor fibrations}, the following strategy was implemented: Consider the germs of \textit{characteristic} vector fields $X^{\lambda}$ at $(\mathbb{R}^{n+1},0)$ associated with \eqref{eq: system 1} given by  
	$$X^{\lambda}:= \sum_{i=1}^{n}a_i^{\lambda}(x,y)\frac{\partial }{\partial x_{i}}+ b^{\lambda}(x,y) \frac{\partial }{\partial y},$$ 
	for each $\lambda=1,\dots,q$, in an analytic coordinate system $(x_1,\dots,x_n,y)$ on an  open set $U\subset \mathbb{R}^{n+1}$.  Suppose that there is a map germ $F_{\eqref{eq: system 1}}:=(f_1,\dots, f_p) :(\mathbb{R}^{n+1},0 ) \to (\mathbb{R}^{p},0)$ admitting only an isolated singularity at the origin, with $1\leq p\leq n+1-q$ where each component $ f_1,f_2\dots, f_p :U \to \mathbb{R} $ is a  \textit{first integral}\footnote{The existence of first integrals has been intensely studied by several mathematicians since the 18th century. Here we do not attempt for further details, but they can be found in the references as  \cite{arnold1,arnold2,Poinca1}.} for \eqref{eq: system 1}, i.e., each $f_j$ is a non-constant function, which is constant on the graph of the solutions of \eqref{eq: system 1} defined on $U$. The map $F_{\eqref{eq: system 1}}$ was called \textit{first integral map germ of system \eqref{eq: system 1}} (for short, \textit{first integral map}). Let $ N_k=f_k^{ -1}(\delta) $ be the \textit{integral manifold} of \eqref{eq: system 1}, i.e., the fibers of the first integral $f_k^{ -1}(\delta)$, for each $k=1,\dots, p $ and a real constant $ \delta $. Denote $M_j:= \bigcap_{k=1}^j N_k$, for any $ j=1,\dots, p$. Since the graph of the solutions $u$ of \eqref{eq: system 1} is contained in the integral manifolds $N_k$, the relations between the topology of Milnor's fibers of $F_{\eqref{eq: system 1}}$ and the topology of the intersection $M_j$ were studied. In other words, the topological classification of the intersections $M_j$ obtained by Milnor Fibration theory was applied to provide a topological description of the graph of the solutions of \eqref{eq: system 1}

	\vspace{0.2cm}
	
	Now, we focus on the more general case where instead of a quasilinear partial differential system, we consider a foliation admitting at least one first integral. Intuitively, a foliation is a decomposition of a manifold into a union of connected, disjoint submanifolds of the same dimension, called leaves, which pile up locally like pages of a book.
	
	\vspace{0.2cm}
	
	More precisely, let $\mathcal{F}$ be a germ of real analytic foliation at $(\mathbb{R}^{n+1},0)$ where $cod(\Sing(\mathcal{F}))\geq 2$. A first integral of $\mathcal{F}$ is a non-constant germ of function $f:(\mathbb{R}^{n+1},0)\to (\mathbb{R},0)$ which is constant on each leaf of $\mathcal{F}$. Then, we consider the map germ $F_{\eqref{eq: system 1}}:=(f_1,\dots, f_p) :(\mathbb{R}^{n+1},0 ) \to (\mathbb{R}^{p},0)$ where $1\leq p\leq n+1-q$ such that each component $ f_1,f_2\dots, f_p :U \to \mathbb{R} $ is a  first integral of $\mathcal{F}$. Furthermore,  $ \dim \Sing F_{\eqref{eq: system 1}} \ge 0 $, i.e., the dimension of the singular set  can be positive. In other words,  $ F_{\eqref{eq: system 1}}  $ may have non-isolated singularities in the fiber over $ \{0\} $.  Additionally, we suppose that first integral map  $F_{\eqref{eq: system 1}}$ preserves the solution of the Laplace equation. This very special class of  maps are known in the literature as \textit{harmonic morphisms}, see e.g., \cite{B1,BW,EL1,F,I,RADG}. 
	
	\vspace{0.2cm}
	
	The authors of \cite{RADG} showed that every polynomial harmonic morphism between Euclidean spaces defines a Milnor fibration in the tube and sphere, and these fibrations are equivalent. They presented a series of new examples of harmonic morphisms and generalized the main result of \cite{PY}. This result guarantees, under certain special hypotheses over a germ of harmonic morphism, that the Milnor sphere fibration is itself a harmonic morphism. This important fact ensures, in particular, the existence of a dictionary between the two theories. Consequently, tools from the theory of harmonic morphisms can be utilized to describe geometric properties of the Milnor fibers of first integral maps.
	
	\vspace{0.2cm}
	
	The main objective of this paper is to show how we can make use of Milnor fibration structures associated with harmonic morphisms between Euclidean spaces to provide a topological and geometric description of the leaves of a foliation admitting at least one first integral.  In particular, we aim to understand the graph of the solutions $u$ of \eqref{eq: system 1}. To achieve this goal, we define in \S2 the \textit{harmonic first integral map} of a foliation and introduce important notions that will be used throughout the text.
	
	\vspace{0.2cm}
	
	
	In \S 3, we will use the relationships between germs of harmonic morphisms and Milnor Fibrations proven in \cite{RADG} to refine the topological description of the  integral manifolds $ N_k$ and the  intersections $ M_j$, for any $ j=1,\dots, p $.

	\vspace{0.2cm}

	In the last section,  we advanced in the understanding of the relationships between germs of harmonic morphisms and Milnor Fibrations,  showing, under the hypothesis of homogeneity of the harmonic first integrals  $F_{\eqref{eq: system 1}}:=(f_1,\dots, f_p) $, that the level set of the  solutions  $u$ of  \eqref{eq: system 1} are contained in spheres whose dimension depends on the number $p$ of first integrals, see Proposition \ref{p1} and Theorem \ref{t2}. Moreover, we will use the relationships between germs of harmonic morphisms and Milnor Fibrations proven in \cite{RADG} and the tools from the theory of harmonic morphisms to describe the  geometric of Milnor fiber associated with a harmonic first integral map, showing that under certain hypotheses, Milnor fibers are minimal manifolds.
	
	\vspace{0.2cm}

	\section{Setup}
	
	\subsection{Some basics facts on foliations}
	Allow us to provide some insight into foliations and integrability. Our first step is the definition of a (nonsingular) foliation.

	\begin{definition}\cite{Camacho1}
		Let $M$ be a $ m $-dimensional manifold. A $p$-dimensional, class $C^r$ $(r\geq 1)$ foliation of $M$ we mean a decomposition of $M$ into a union of disjoint connected subset $\{L_{\alpha}\}_{\alpha \in A}$, called the leaves of the foliation, with the following property: Every point in $M$ has a neighborhood $U$ and a system of local, class $C^r$ coordinates $x=(x^1,\dots,x^m):U\to \mathbb{R}^m$ such that for each leaf $L_{\alpha}$, the components of $U\cap L_{\alpha}$ are described by the equations $x^{p+1} =constant, \dots, x^m=constant$.  	
	\end{definition}

	A distribution of $k$-planes on a manifold $M$ is defined as a mapping $D$ that assigns to each point $x \in M$ a vector subspace of dimension $k$ in $T_xM$. In other words, for every $q \in M$, there exist $k$ vector fields $X_1, \ldots, X_k$  (non-vanishing) on $M$ such that ${X_1(x), \ldots, X_k(x)}$ forms a basis for $D(x)$. If, for any two vector fields $X$ and $Y$ on $U$ such that $X(x)$ and $Y(x)$ belong to $D(x)$ for each $x \in U$, it follows that $X,Y$ also belongs to $D(x)$, then the distribution $D$ is said to be involutive. The Frobenius Theorem (for more details, see for example \cite{Camacho1}) guarantees that an involutive distribution $D$ defines a foliation on $M$ of dimension $k$.
	
	\vspace{0.2cm}
	
	Now, let $\omega_{1}, \ldots, \omega_{q} \in \Omega^{1}(M)$ be linearly independent $1$-forms (non-vanishing) on $M$. This gives rise to a corresponding distribution $D$ of $(n-q)$-dimensional planes, defined as
	\[
	D(x) = \{v\in T_pM, \omega_j(x)\cdot v =0,j=1,\dots,q\},
	\]
	for each $x\in M$. 
	
	\begin{definition}
		We say that the system $\omega_{1},\ldots, \omega_{q}$ is integrable if it satisfies
		\begin{equation*}
			d\omega_{j} \wedge \omega_1 \wedge \omega_2 \wedge \dots \widehat{\omega_j}\dots \wedge\omega_q = 0 \quad \text{for all } j=1,\ldots, q \quad \text{on } M.
		\end{equation*}
	\end{definition}
	
	It follows from the Frobenius Theorem that the system of $1$-forms $\omega_1,\dots,\omega_q$ is integrable if and only if the corresponding distribution $D$ defines a foliation on $M$ of dimension $n-q$.

	\vspace{0.2cm}

	A singular foliation can be defined as a pair $\mathcal{F} = (\mathcal{F}^*, S)$, where $S\subsetneq M$ is a
	subset and $\mathcal{F}^*$ is a nonsingular holomorphic foliation of $M\setminus S$. We refer to $S$ as the singular set of
	$\mathcal{F}$ and denote it as $S = \Sing (\mathcal{F})$. The leaves of $\mathcal{F}$  are the leaves of $\mathcal{F}^*$  on $M\setminus \Sing (\mathcal{F})$. In what follows, we will consider  $\codim(\Sing(\mathcal{F}))\geq2$.

	\vspace{0.2cm}

	Let us point out that $k$ vector fields $X_1, \dots, X_k$ can define a singular foliation of dimension $k$, where $\Sing(\mathcal{F}) =\{p \in M : X_1\wedge \dots \wedge X_k(p) =0 \}$. Moreover, the $1$-forms $\omega_1,\dots, \omega_{m-k}$ can define a singular foliation of dimension $k$ by $m-k$ $1$-forms such that $\Sing(\mathcal{F}) = d\omega_i\wedge \omega_j \equiv 0$.

	\subsection{First integral maps of real analytic foliations}\label{section: Pre}

	\begin{definition}\label{d1} 
		Let $\mathcal{F}$ be a germ of (possible singular) real analytic foliation of dimension $q<n+1$ at $(\mathbb{R}^{n+1},0)$. A \textit{first integral} of $\mathcal{F}$ is a germ of non-constant function $f:(\mathbb{R}^{n+1},0)\to (\mathbb{R},0)$ which is constant on each leaf of $\mathcal{F}$. A \textit{first integral map germ} (for short, \textit{first integral map}) of the foliation $\mathcal{F}$ is a map germ $F:=(f_1,f_2\dots,f_p):(\mathbb{R}^{n+1},0)\to (\mathbb{R}^{p},0)$ such that $f_1,\dots,f_p$ are  functionally independent\footnote{Let us remind that a set of functions $  f_1,\dots,f_k   $  are \textit{functionally independent} if $df_1\wedge \dots \wedge df_k \not\equiv 0$.} first integrals of $\mathcal{F}$. 
	\end{definition}
	
	\vspace{0.2cm}
	
	We say that $f_1,\dots,f_{p}$ is a \textit{complete system of first integrals} of $\mathcal{F}$ if $f_1,\dots,f_p  $  are functionally independent first integrals of $\mathcal{F}$ where  $p=n+1-q$. In this case we say that $\mathcal{F}$ is a \textit{complete intersection}.
	
	\vspace{0.2cm}
	
	The singular set of $ F$, is given by $\Sing F :=\{x\in U: \rank( DF(x)) <p\}$, where $DF(x)$ denotes the Jacobian matrix of $F$ at $x$, and $U \subset \mathbb{R}^m$ is an open set with $0\in U$. 
	
	
	
	\vspace{0.2cm}

	We say that $F$ has an \textit{isolated singularity at the origin} when $\Sing F = \{0\}$. This condition means that there exists a neighborhood $ U $ of the origin in $ \mathbb{R}^{n+1} $ such that the jacobian matrix $DF(x)$ has rank $ p $ at all points $ x \in U $ other than the point $ x=0 $. In this work, we will focus on the least restrictive condition, where $\rm{dim} \, \Sing F \ge 0$, with it is still included in $ V_F:=F^{-1}(0) $. In this last case, we denote $ \Disc F := F(\Sing F) =\{0\}$, and we say that $ F $ has an \textit{isolated critical value.}
	
	\vspace{0.2cm}

	\begin{example}\cite{RRE2} \label{ex1} Let $X_1,X_2$ be two germs of vector fields at $(\mathbb{R}^4,0)$ determined by 
		$$X_{1}= 	2 x_3 x_4\frac{\partial}{\partial x_2} - 3x_4(x_1^2 +x_2^2)\frac{\partial}{\partial x_3} -   3x_3(x_1^2 +x_2^2)\frac{\partial}{\partial x_4}$$ 
		and 
		$$X_2=	x_4\frac{\partial }{\partial x_2} + x_4^3 \frac{\partial }{\partial x_3}  -\left(\frac{3}{2}(x_1^2+x_2^2)+x_3x_4\right)\frac{\partial }{\partial x_4}.$$ 
		The above vector fields determine a codimension two foliation $\mathcal{F}$ at $(\mathbb{R}^4,0)$. This foliation also can be determined by the integrable $2$-form $\omega$ given by
		\[
		\omega = (3x_1^2+3x_2^2)dx_1\wedge dx_2 + 2x_3 dx_1\wedge dx_3 +2x_4 dx_1\wedge dx_4.
		\] Now, define the functions $f_1(x_1,x_2,x_3,x_4)=x_1$ and $f_2(x_1,x_2,x_3,x_4)=3x_1^2x_2+x_2^3+x_3^2+x_4^2$. Through a straightforward computation we obtain $df_i(X^\lambda)=0$ for every $i,\lambda \in\{1,2\}$ and $df_1\wedge df_2 \neq 0$. Furthermore, we have $\omega = df_1 \wedge df_2$. Since $ p=2,q=2,n=3$ then $p= n+1-q $ and $F=(f_1,f_2)$ is first integral map such that the components $f_1,f_2$ are complete system of first integral for $\mathcal{F}$. Note that $\Sing(F)= \{0\}$.
	\end{example}

	\vspace{0.2cm}
	\begin{example}\cite{RRE2}  Let $\mathcal{F}$ be a germ of one-dimensional foliation at $(\mathbb{R}^3,0)$ determined by the vector field 
		\[
		X= (3x_2x_3^2+x_3)\frac{\partial}{\partial x_1} + x_1x_4(3x_4-2)\frac{\partial}{\partial x_2} - (2x_2+1)x_1\frac{\partial }{\partial x_3}.
		\]
		The map $F=(f_1,f_2):(\mathbb{R}^3,0)\to(\mathbb{R}^2,0)$ where $f_1(x_1,x_2,x_3) = x_1^2-x_2^2+x_3^2$ and $f_2(x_1,x_2,x_3) = x_1^2+x_2+x_3^3$ is a first integral map of the $\mathcal{F}$ with $\dim \Sing(F) > 0$. 	
	\end{example}
	
	\vspace{0.2cm}
	
	\begin{example} Consider the real analytic germs of $1$-forms $\omega_1,\omega_2$ at $(\mathbb{R}^8,0)$ given by
		\begin{equation*}\label{eq: eq2}
			\begin{split}
				\omega_1 & = (3x_1^2 x_7 + 6x_1x_2x_8-3x_2^2x_7)dx_1 + (3x_8x_1^2-6x_1x_2x_7-3x_2^2x_8)dx_2 \\
				&-2x_3dx_3 +2x_3dx_4 -3(x_6^2-x_5^2)dx_5  - 6x_5x_6dx_6 \\&+ (x_1^3-3x_1x_2^2)dx_7 + (3x_1^2x_2-x_2^3)dx_8 
			\end{split}
		\end{equation*}
		and 
		\begin{equation*}
			\begin{split}
				\omega_2 & = (3x_1^2x_8-6x_1x_2x_7-3x_2^2x_8)dx_1-(3x_7x_1^2+6x_1x_2x_8-3x_2^2x_7)dx_2 \\ &+ 2x_4dx_3 + 2x_3dx_4 +2x_5x_6dx_5 + (x_5^2-3x_6^2)dx_6 \\ &- (3x_1^2x_2-x_2^2)dx_7 + (x_1^3-3x_1x_2^2)dx_8,
			\end{split}
		\end{equation*}
		Since $\omega_j\wedge d\omega_j \equiv 0, \forall j\in \{1,2\}$, then $\omega_1, \omega_2$ define two condimension one (real analytic) foliations $\mathcal{F}_1,\mathcal{F}_2$ at $(\mathbb{R}^8,0)$, respectively. Now, define the functions $f_1,f_2:(\mathbb{R}^8,0)\to(\mathbb{R},0)$ given by $$f_1(x_1,\dots,x_8)=x_7x_1^3-3x_7x_1x_2^2+3x_8x_1^2x_2-x_8x_2^3+x_5^3-3x_5x_6^2+x_4^2-x_3^2$$ and $$f_2(x_1,\dots,x_8)=-3x_7x_1^2x_2+x_7x_2^3+x_8x_1^3-3x_8x_1x_2^2+x_5^2x_6-x_6^3+2x_4x_3.$$  
		Since $\omega_j = df_j$ then $f_j$ is a first integral for $\mathcal{F}_j$. Furthermore, since $df_1\wedge df_2\not\equiv 0$ we can define the codimension two foliation $\mathcal{F} = \mathcal{F}_1\cap \mathcal{F}_2$ at $(\mathbb{R}^8,0)$. The first integral map $ F=(f_1,f_2) $ satisfies $ \dim \,\Sing F = 2 $ and $ \Disc F= \{0\} $.  	
	\end{example}

	\subsection{Harmonicity and first integral maps}
	
	Let $M$ be a Riemannian manifold and we consider $\Omega^{k}(M)$ the space of $k$-forms on $M$ where 
	$$d^k\colon \Omega^{k}(M)\rightarrow \Omega^{k+1}(M),$$
	is exterior fifferentiation. The inner product structures on $\Omega^k(M)$ allow us to define the adjoint 
	$\delta^k\colon \Omega^{k+1}(M)\rightarrow \Omega^k(M)$, to $d^k$ via the formula $\left(\delta^k\omega_{1},\omega_{2}\right)=\left(\omega_{1},d^k\omega_{2}\right)$.
	From the definition of the Hodge star operator, we see that $\delta$ can be computed from $d$ as follows 
	$\delta^k=(-1)^{(n-k)(k+1)}\star d^{n-k-1}\star$.
	Consequently,  the Laplacian on forms\footnote{also called the Hodge Laplacian}  $\Delta\colon \Omega^{K}(M)\rightarrow \Omega^{K}(M)$, is given by
	$\Delta\omega=(d\delta+\delta d)\omega$.
	
	\vspace{0.2cm}
	
	\begin{definition}
		We say that a $1$-form $\omega$ is \textit{harmonic} if $\Delta \omega=0$, i.e., $\omega$ is closed and co-closed.
	\end{definition}

	\vspace{0.2cm}
	
	Let $\omega$ be a germ of $1$-form at $(\mathbb{R}^{n+1},0)$. If $\omega$ is harmonic then it follows from Poincaré-Lemma that $\omega$ is exact. Therefore, there is a harmonic function $f$ such that $\omega = df$. Now, consider the set of germs of harmonic $1$-forms $\omega_1,\dots,\omega_q$ at $(\mathbb{R}^{n+1},0)$. 
	It it follows from the above argument that there are harmonic functions $f_1,\dots,f_q$ such that $\omega_j = df_j$ for every $j\in\{1,\dots,q\}$. Let $\omega$ be the $q$-form given by $\omega = \omega_1\wedge \dots \wedge \omega_q$ where each $\omega_j$ is a harmonic $1$-form.  Then, $\omega$ determines a codimension $q$ foliation $\mathcal{F}$ at $(\mathbb{R}^{n+1},0)$. Note that $ F=(f_1,f_2\dots, f_p) :(\mathbb{R}^{n+1},0) \to (\mathbb{R}^{p},0)$ is a first integral map of $\mathcal{F}$ satisfying the Laplace equation $ \sum_{i=1}^{n+1}\dfrac{\partial^2 f_\alpha}{\partial x_{i}^{2}} = 0 $, for any $ \alpha=1,\ldots,p. $ Note that this is the case of all examples in Section \ref{section: Pre}. 
	
	\vspace{0.2cm}
	
	With the discussion presented above as our motivation, we now introduce the following definition, which describes the type of foliation we are dealing with. Specifically, we are interested in foliations that possess at least one first integral that is a harmonic function.
	
	\vspace{0.2cm}
	
	\begin{definition}\label{dhfim} Let $\mathcal{F}$ be a $q$-dimensional germ of real analytic foliation at $(\mathbb{R}^{n+1},0)$. Suppose that $\mathcal{F}$ admits $f_1,\dots,f_p$ functionally independent first integrals, where $p\geq 1$. We say that $$ F=(f_1,f_2\dots, f_p) :(\mathbb{R}^{n+1},0) \to (\mathbb{R}^{p},0), $$ is a  \textit{harmonic first integral map}, if $ F $ is a first integral map for $ \mathcal{F} $ and  satisfies:
		\begin{enumerate}
			\item [(1)] $ \sum_{i=1}^{n+1}\dfrac{\partial^2 f_\alpha}{\partial x_{i}^{2}} = 0 $, for any $ \alpha=1,\ldots,p $  \medskip;
			
			\item[(2)] $ \sum_{i=1}^{n+1}\dfrac{\partial f_\alpha}{\partial x_{i}}\dfrac{\partial f_\beta}{\partial x_{i}} = \lambda^2(x)\delta_{\alpha\beta}$, for any $ \alpha, \beta = 1, \ldots, p,$
		\end{enumerate}
		where $\delta_{\alpha\beta}$ is the delta of Kronecker, and $ \lambda: (\bR^{n+1},0) \to (\bR,0) $ is a germ of continuous function. 
	\end{definition}
	
	\begin{remark}
		It follows from condition (2) that  the critical locus of a harmonic  first integral map $F$ is  given by $ \Sing F=\left\lbrace x\in \bR^{n+1}\,|\,\|\nabla f_j(x)\| = 0, \textrm{for any fixed } j \right\rbrace.$    Therefore, if $ x\in \Sing F$ then $ \nabla f_j(x)=0 $ for any $ j=1,\ldots,p $, i.e., $DF(x) \equiv 0$ .
		
	\end{remark}	
	
	\vspace{0.2cm}
	
	Let us point out that after  \cite{F} and \cite{I}, the conditions $ (1) $ and $ (2) $ in the Definition  \ref{dhfim} state that a first integral map germ $ F=(f_1,f_2\dots, f_p) :(\mathbb{R}^{n+1},0) \to (\mathbb{R}^{p},0), $ is a {\it harmonic morphism} i.e., a map between Riemannian manifolds that preserve the solutions of the Laplace equation in the sense that for any function $ f: W \to \bR $ which is harmonic on an open set $ W  \subset \bR^p$, with $ F^{-1} (W)$ non-empty,  the composition $$ f\circ F: F^{-1} (W) \to \bR $$ is harmonic on $ F^{-1} (W) \subset \bR^{n+1} $, i.e.,  $ F$  pullback germs of harmonic functions on $ \bR^p $ to germs of harmonic functions on $\bR^{n+1} $. In other words, the conditions $ (1) $ and $ (2) $  are equivalent to saying that $ F$ is  both a {\it harmonic map} and {\it horizontally weakly conformal map}, respectively,  where harmonic  map is a  critical point of the energy functional $$ E(F,\Omega)=\dfrac{1}{2} \int_{\Omega} |\textrm{d}F|^2$$ 
	over any compact set $ \Omega \subset \bR^{n+1} $ and the property horizontally weakly conformal, means that at the point where $ \dd F(x) \not = 0 $, $\dd F(x)$ preserves horizontal angles. For more details, see e.g., \cite{B1,BW,EL1,F,I}. 
	
	%
	


	\vspace{0.2cm}
	
	The next result is a particular case of \cite[Corollary 2.3]{ABB}, proved by Ababou, R., P. Baird and J. Brossard. It shows that if the first integrals of a foliation are polynomial functions, and their gradients are mutually orthogonal and of equal length $ \lambda(x) $ at each point $ x\in \bR^{n+1} $, then $ F$ is a harmonic first integral map.  
	
	\vspace{0.2cm}
	
	\begin{theorem}\cite[Corollary 2.3]{ABB}\label{phm}
		Consider a first integral map  $ F:=(f_1,f_2\dots, f_p) :(\mathbb{R}^{n+1},0) \to (\mathbb{R}^{p},0) $ of foliation $\mathcal{F}$. If $ F $ is a horizontally weakly conformal map and  $f_1,\dots,f_p$  are polynomial functions, then it is a harmonic first integral map.
	\end{theorem}
	
	\vspace{0.2cm}
	
	\begin{example}\label{ex: dimSing positive}
		Let $\mathcal{F}$ be the germ of codimension two foliation at $(\mathbb{R}^4,0)$ determined by the $2$-form $\omega = \omega_1 \wedge \omega_2$ where
		\begin{equation*}
			\begin{split}
				\omega_1 &=( -2w^3y-6w^2xz+6wyz^2+2xz^3 )dx +(-2w^3x+6w^2yz+6wxz^2-2yz^3 )dy \\
				& +(-3w^2x^2+3w^2y^2+12wxyz+3x^2z^2-3y^2z^2)dz \\
				& + (-6w^2xy-6wx^2z+6wy^2z+6xyz^2 ) dw
			\end{split}
		\end{equation*}
		and
		\begin{equation*}
			\begin{split}
				\omega_2 &=(2w^3x-6w^2yz-6wxz^2+2yz^3)dx + (-2w^3y-6w^2xz+6wyz^2+2xz^3 )dy \\
				& +(-6w^2xy-6wx^2z+6wy^2z+6xyz^2)dz \\
				& +(3w^2x^2-3w^2y^2-12wxyz-3x^2z^2+3y^2z^2)dw.
			\end{split}
		\end{equation*}
		Now, define the map $ F=(f_1,f_2):(\mathbb{R}^4,0) \to (\mathbb{R}^2,0) $, where
		$$ \left\lbrace \begin{array}{ccl}
			f_1(x,y,z,w)&= -2w^3xy-3w^2x^2z+3w^2y^2z+6wxyz^2+x^2z^3-y^2z^3 &   \medskip\\
			f_2(x,y,z,w)&= w^3x^2-w^3y^2-6w^2xyz-3wx^2z^2+3wy^2z^2+2xyz^3. &  \medskip\\

		\end{array}\right. $$
		One has that $df_j = \omega_j, j\in \{1,2\}$. Then, $F$ is a first integral map for $\mathcal{F}$. Furthermore, $ \left\langle \nabla f_1, \nabla f_2 \right\rangle =0 $ and $ \|\nabla f_1 \|^2 =(x^2+y^2)(w^2+z^2)^2(4w^2+9x^2+9y^2+4z^2)= \|\nabla f_2 \|^2  $. Thus, $ F $ is a horizontally weakly conformal map. Since $ F $ is polynomial, it follows from Theorem \ref{phm}  that $ F$ is a harmonic first integral map for $\mathcal{F}$. 

	\end{example}

	\vspace{0.2cm}
	
	\begin{example}\label{ex4-3}
		Consider the germ of $2$-form $\omega$ at $(\mathbb{R}^4,0)$ given by
		\begin{equation*}
			\begin{split}
				\omega  &= (x_1x_4-x_2x_3)dx_1 \wedge dx_2 +  (x_1^2+x_3^2) dx_1\wedge dx_3 +(x_1x_2+x_3x_4)dx_1\wedge dx_4 \\ &+ (x_1x_2+x_3x_4)dx_2\wedge dx_3 + (x_2^2+x_4^2)dx_2\wedge dx_4 + ( x_1x_4-x_2x_3)dx_3\wedge dx_4.
			\end{split}
		\end{equation*}
		We have that $\omega$ determines a germ of foliation at $(\mathbb{R}^4,0)$ admitting a first integral map $F=(f_1,f_2):(\mathbb{R}^4,0) \to (\mathbb{R}^2,0) $, where
		$$ \left\lbrace \begin{array}{ccl}
			f_1(x,y,z,w)&= x^2+y^2-z^2-w^2 &   \medskip\\
			f_2(x,y,z,w)&= 2(xz+yw) &  \medskip\\
		\end{array}\right. $$
		Furthermore, $ \left\langle \nabla f_1, \nabla f_2 \right\rangle =0 $ and $ \|\nabla f_1 \|^2 = \|\nabla f_2 \|^2  $. Thus, $F$ is a horizontally weakly conformal map. Since $ F $ is polynomial, it follows from Theorem \ref{phm}  that $ F $ is a harmonic morphism. This map is known as the \textit{Hopf polynomial map}.

	\end{example}

	\section{The Milnor fibrations for  harmonic first integral maps}
	
	\vspace{0.2cm}
	
	In this section, we will show how the topology of Milnor fibers associated with a harmonic first integral map can be used to give a topological description of the leaves of the foliation. First, let's recall some basic facts about the existence of Milnor fibrations. The existence of such fibrations as well its equivalence problem has been studied by several mathematicians since 80' decade. Here we do not attempt for further details, but they can be found in the references as \cite{Mi,Le,Ma,PT,ART,dST0,dST1,A,ACT,CSS3,CSS1,PS0,PS2,PS3,RSV,RA,S1}.

	\vspace{0.2cm}


	Let $F:(\bR^{n+1},0) \to (\bR^p,0),$ $m\geq p\geq 2,$  be a real analytic map germ and let $F: U \to \bR^{p},$ $F(0)=0$ be an analytic representative of the germ where $U \subseteq \mathbb{R}^{n+1}$ is an open set and $F(x)=(f_1(x),f_2(x), \ldots,f_p(x)).$  The singular and discriminant sets of $ F $, denoted by $\Sing(F)$ and $\Disc(F)$, respectively, are defined as was done for the first integral map, see Definition \ref{d1}. We remind  that $F$ has an isolated critical value at origin if $\Disc(F)=\{0\}.$
	
	\vspace{0.3cm}
	
	\begin{definition}\cite{ART,Ma,Mi,Le,PT}\label{def: Milnor's tube} We say that a real analytic map germ $F:(\bR^{n+1},0) \to (\bR^p,0),$ admits a Milnor tube fibration if there exist $\epsilon_{0}>0$ and $\eta>0,$ such that for all $0< \epsilon \leq \epsilon_{0}$ in the  ball ${B}_{\epsilon}$ the restriction map 	
		$$F_{|}:{B}^{n+1}_{\epsilon}\cap F^{-1}(S_{\eta}^{p-1})\to S_{\eta}^{p-1}$$
		
		\noindent is a smooth locally trivial  fibration, for all $0<\eta \ll \epsilon.$ Moreover, the diffeomorphism type does not depends on the choice of $\epsilon$ and $\eta$.
		
	\end{definition}
	
	\vspace{0.3cm}

	The existence of the Milnor tube fibration is well known and in order to state them we need to introduce some notations and definitions.
	
	\vspace{0.3cm}
	
	Let $U \subset \mathbb{R}^m$ be an open set,  $0\in U$,  and let $\rho:U \to \mathbb{R}_{\ge 0}$ be  the function square of the Euclidean distance for the origin, i.e., $\rho(x)=\|x\|^{2}$. We consider the following definition:
	
	\vspace{0.2cm}
	
	\begin{definition}\label{d:ms}
		Let $F:(\mathbb{R}^{n+1}, 0) \to (\mathbb{R}^p,0)$ be a non-constant analytic map germ. The set germ at the origin:
		\[M(F):=\Sing(F, \rho) \]
		is called the set of \textit{$\rho$-nonregular points} of $F$, or the \emph{Milnor set of $F$}.
	\end{definition}

	\vspace{0.2cm}
	
	The following condition was first used in \cite{Ma} and later in \cite{ACT} to ensure the existence of the Milnor tube fibrations.
	\begin{equation}\label{eq:main}
		\overline{M(F)\m F^{-1}(0)}\cap F^{-1}(0) =\{0\}.
	\end{equation}
	
	\vspace{0.2cm}

	The result about the existence of a Milnor tube fibration can be state for a harmonic first integral map of a foliation as follows: 
	
	\vspace{0.2cm}

	\begin{theorem}\label{tf1} If  $ F=(f_1,f_2\dots, f_p) :(\mathbb{R}^{n+1},0) \to (\mathbb{R}^{p},0) $, $p \geq 2,$  is  a harmonic  first integral map of a foliation $\mathcal{F}$, then  admits a Milnor tube fibration
		\begin{equation}\label{tfeq1}
			{F}_|:{B}^{n+1}_{\epsilon}\cap F^{-1}(S_{\eta}^{p-1})\to S_{\eta}^{p-1}
		\end{equation}
		
	\end{theorem}
	
	\begin{proof}
		The result follows from Definition \ref{dhfim} and \cite[Theorem 3.1]{RADG}, once every harmonic morphism satisfies the condition \eqref{eq:main} and has an isolated critical value  at the origin. 
	\end{proof}
	
	\vspace{0.2cm}
	
	\begin{remark} Let us point out that if $ f:(\mathbb{R}^{n+1},0) \to (\mathbb{R},0) $ is a analytic first integral of a foliation, then $ \Disc f = \{0\} $ and $f$ satisfies the condition \eqref{eq:main}. For more details, see \cite{DA}.
		
			
			
		
	\end{remark}
	
	
	Next, let us remind the definition of Milnor sphere fibrations.
	
	\begin{definition} \label{def: Milnor's sphere} We say that a real analytic map germ $F:(\bR^{n+1},0) \to (\bR^p,0),$ admits a  Milnor sphere fibration if there exists $\epsilon_{0}>0$ such that for all $0< \epsilon \leq \epsilon_{0}$ the  projection map 
		
		\begin{equation}\label{smil}
			{\Psi_F} := \dfrac{F}{\|F\|}:S_{\epsilon}^{n}\m K_{\epsilon}\to S^{p-1}
		\end{equation}
		
		\noindent
		is a locally trivial smooth fibration, for all $0<\epsilon \ll \epsilon_{0}.$ The closed set $K_{\epsilon}:=V_{F}\cap S_{\epsilon}^{n}$ is called the \textit{link} of the singularity of $F$.
		
	\end{definition}
	
	Given $F :U\to \mathbb{R}^{p}$ as above consider the map projection:
	\begin{equation}\label{eq:spherefib}
		{\Psi_F} := \frac{F}{\| F\|}: U \m V_F \to  S^{p-1}.
	\end{equation}
	
	We define the Milnor set $M({\Psi_F}):=\Sing ({\Psi_F}, \rho)$ of the map \eqref{eq:spherefib}, i.e. the germ at the origin of the $\rho$-nonregular points of ${\Psi_F}$. We say that \textit{$ \Psi_F$ is $\rho$-regular} if $M ({\Psi_F}) =\emptyset.$
	
	\vspace{0.2cm}

	The result about the existence of a Milnor sphere fibration can be state for a harmonic first integral map of a foliation as follows: 
	
	\vspace{0.2cm}

	\begin{theorem}\label{tf3} Let  $ F=(f_1,f_2\dots, f_p) :(\mathbb{R}^{n+1},0) \to (\mathbb{R}^{p},0) $, $p \geq 2,$  be  a harmonic  first integral map of a foliation $ \mathcal{F} $. Then $F$ admits the Milnor sphere fibration
		\begin{equation}\label{tfeq2}
			\Psi_{F} :S_{\epsilon}^{n}\m K_{\epsilon}\to S^{p-1}.
		\end{equation}
		
	\end{theorem}
	
	\begin{proof}
		It  follows from Definition \ref{dhfim} and \cite[Theorem 3.2]{RADG}, once the projection  $\Psi_{F}$  is $ \rho $-regular. 	
	\end{proof}

	\vspace{0.2cm}
	
	\begin{example}\label{ex5-2}
		Let $\mathcal{F}$ be the codimension two germ of foliation at $(\mathbb{R}^4,0)$ determined by $\omega=\omega_1\wedge \omega_2$ where
		\[
		\omega_1 = (a^2-b^2)dx+(-2ab)dy +(2ax-2by)da +(-2ay-2bx)db
		\] 
		and
		\[
		\omega_2 = 2abdx + (a^2-b^2)dy+ (2ay-2bx)da+   (2ax-2by)db.
		\]
		We have that $\mathcal{F}$ 	admits a first integral map $ F=(f_1,f_2):(\mathbb{R}^4,0) \to (\mathbb{R}^2,0) $, where
		$$ \left\lbrace \begin{array}{ccl}
			f_1(x,y,a,b)&=a^2x-2aby-b^2x  &   \medskip\\
			f_2(x,y,a,b)&= a^2y+2abx-b^2y. &  \medskip\\
		\end{array}\right. $$  Furthermore, 	$ \left\langle \nabla f_1, \nabla f_2 \right\rangle =0 $ and $ \|\nabla f_1 \|^2 =(a^2+b^2)(a^2+b^2+4x^2+4y^2)= \|\nabla f_2 \|^2  $. Thus, $ F$ is a horizontally weakly conformal map. Since $ F$ is polynomial, it follows from Theorem \ref{phm}  that $ F$ is a harmonic morphism. Now, it is follows from  Theorems \ref{tf1} and \ref{tf3},  	the harmonic first integral map $ F$ admits a Milnor tube and sphere  	fibrations.
		

	\end{example}

	\subsection{Topology of the leaves}

	Let $\mathcal{F}$ be a foliation admitting a first integral $f$. Then, outside \Sing($\mathcal{F})$ we have that $f$ is constant on each leaf $L$ of $\mathcal{F}$. In particular, $L\subset f^{-1}(c)$, for some constant $c\in \mathbb{R}$. In other words,  each leaf of a foliation lies in the level sets of their first integrals, i.e., in a fiber of $ f $. However, in the neighborhood of the singular set of the foliation, the behavior of the leaves can be much more complicated and the above relation may not exist.
	
	\vspace{0.2cm}
	
	Our main objective is to provide more information about the topology and geometry of the leaves that are in proximity to the singular set of a given foliation. To achieve this, we employ the Milnor fibration associated with the harmonic first integral maps of the foliation. 
	

		\vspace{0.2cm}
		
		Suppose that there is a harmonic first integral map $F=(f_1,f_2\dots, f_p) :(\mathbb{R}^{n+1},0) \to (\mathbb{R}^{p},0) $, $p \geq 2,$  of $\mathcal{F}$. The Milnor fiber for \eqref{tfeq1} is given by $ M_{F}:= F^{-1}(y) $, where $ y=(y_1, \dots, y_p) \in S_{\eta}^{p-1}  $. It is well known that the Milnor fiber $ M_F$ is diffeomorphic to the manifold
		$ {B}^{n+1}_{\epsilon} \cap M_p $, where  $M_p:=\bigcap_{k=1}^p N_k$, i.e.,  the intersection of the integrals manifold  $N_k = f_k^{-1}(\delta_k) $, for any $ k=1,\dots, p $,  and $ \delta_k $  arbitrary  small enough real numbers, not all equal to zero. 
		
		\vspace{0.2cm}
		
		The Theorem 3.6 of \cite{RADG} shows that Milnor fibration in the tube and sphere are equivalent, where \textit{equivalence}, means that there is a diffeomorphism $ h: S_{\epsilon}^{n}\m K_{\epsilon} \to \left( {B}^{n+1}_{\epsilon}\cap F^{-1}(S_{\eta}^{p-1})\right)  $ such that $ \Psi_{F} ={F}_{|} \circ h   $. Consequently, $ h $ induces a diffeomorphism $ h_y:\Psi_{F}^{-1}(y) \to  M_F$,  for any $ y=(y_1, \dots, y_p) \in S_{\eta}^{p-1}  $, i.e., the Milnor fiber on the tube and on the sphere associated to a harmonic first integral map, are diffeomorphic.
		
		\vspace{0.2cm}
		
		As it is possible to verify directly, if $  F $ is a harmonic first integral map  then the  map  $ F^j:=(f_1,f_2\dots, f_j) :(\mathbb{R}^{n+1},0) \to (\mathbb{R}^{j},0) $, with $ 2\le j < p $ is a harmonic first integral map, as well. Consequently,  $ F^j$ has a Milnor tube and sphere fibrations and again, the Milnor fibers $ \Psi_{F^j}^{-1}(y) $ and $ M_{F^j} $ are diffeomorphic to the manifold
		$ {B}^{n+1}_{\epsilon} \cap M_j $, with $M_j:=\bigcap_{k=1}^j N_k$.

		\vspace{0.2cm}
		
		\subsubsection{The isolated singular case}
		\vspace{0.2cm}
		
		Let  us recalling a framework which was approach in \cite{RRE2}. Consider a  system  of quasi-linear partial differential equations of first-order given by
		\begin{equation}\label{eq: system prop}
			\sum_{i=1}^{n}a_i^{\lambda}(x,u)\frac{\partial u}{\partial x_{i}}= b^{\lambda}(x,u) \tag{$ \ast $},
		\end{equation}
		where $u$ is one unknown function in $n$ independent variables $x=(x_1,\dots,x_n)$, $a_i^{\lambda}$ and $b_i^{\lambda}$  nonzero analytic functions defined on an open set $\Omega \subset \mathbb{R}^{n+1}=\{(x_1,\dots,x_n,u)\}$, for $ i=1\dots,n$ and $\lambda =1,\dots,q$, with $1\leq q\leq n$. Sometimes, for the sake of simplicity we will denote $x_{n+1}:=u$ and $a_{n+1}^\lambda := b^{\lambda}$. Then we consider the \textit{characteristic} vector fields associated with the quasilinear system given by $$X^{\lambda}:= \sum_{i=1}^{n}a_i^{\lambda}(x,u)\frac{\partial }{\partial x_{i}}+ b^{\lambda}(x,u) \frac{\partial }{\partial u},$$ for each $\lambda=1,\dots,q$. Note that $X^\lambda$ define a real analytic foliation of dimension $q$, since $[X^{\lambda_1},X^{\lambda_2}]\in D$ for every $\lambda_1,\lambda_2$ and $X^{\lambda_1}\wedge\dots \wedge X^{\lambda_q}\not\equiv 0$. Then, a first integral of the system \eqref{eq: system prop}  is a non-constant function $f:\Omega\to \mathbb{R}$ such that $f$ is constant on the general solutions of \eqref{eq: system prop}  defined on $\Omega$. Then considering the characteristic vector fields $X^\lambda$  one has $df(X)\equiv 0$, i.e., on $\Omega$, one has $$\sum_{i=1}^{n}a_i^{\lambda}(x,u)\frac{\partial f(x,u)}{\partial x_{i}}+ b^{\lambda}(x,u) \frac{\partial f(x,u)}{\partial u}=0.$$
		Let $f_1,\dots,f_p$  be first integrals functionally independent  for the system \eqref{eq: system prop} for some $p\geq 1$. The map germ $$ F_{\eqref{eq: system prop}}:=(f_1,f_2\dots, f_p) :(\mathbb{R}^{n+1},0) \to (\mathbb{R}^{p},0), $$  is called  the  \textit{first integral map germ of system \eqref{eq: system prop}}, for short, \textit{first integral map}. In this case, $f_1,\dots,f_{p}$ is a \textit{complete system of first integrals} for the system \eqref{eq: system prop}, if $  f_1,\dots,f_p  $  are functionally independent first integrals and  $p=n+1-q$. 
		Since the level sets of the solutions $u$ of  \eqref{eq: system prop}  are contained in integral manifolds (see for instance Chapter 1 in \cite{McOwen}), following the strategy of the authors of \cite{RRE2}, if we have a sufficient number of first integrals $  f_1,f_2\dots, f_p $ and we know the integral manifolds $ N_k$, one can implicitly describe the solutions of the system, provided we are able to describe locally, the topology of the intersection $ M_j$, for any $ j=1,\dots, p $. Therefore, from the above discussion, it is exactly at this point that we can make use of the existence of local fibration structures in neighborhood of  singularities of harmonic first integral maps, once this  allows the description of the respective Milnor fibers 
		and consequently, the description of the intersections $ M_p $ and $ M_j $.
		
		Also in \cite{RRE2}, the authors applied the results of \cite{DDA,Iodim,SHMA} and \cite{CL} and stated the following results for first integral maps for the system \eqref{eq: system prop} with an isolated singularity at the origin, which, in the context of harmonic first integral map, reads:
		
		\vspace{0.2cm}
		
		\begin{theorem}\cite{DDA,RRE2}\label{t1}
			Let  $ F_{\eqref{eq: system prop}}=(f_1,f_2\dots, f_p) :(\mathbb{R}^{n+1},0) \to (\mathbb{R}^{p},0) $ be  a harmonic first integral map for the system \eqref{eq: system prop} , 
			which has an isolated singularity at the origin. The following holds true:
			\begin{enumerate}
				\item [(i)] If $ n+1 $ is even, then $ \chi(M_p) = 1-\deg_0 \nabla f_1 $.  Moreover, one has $ \deg_0 \nabla f_1=\deg_0 \nabla f_2 = \ldots = \deg_0 \nabla f_p,  $ where  $ \deg_0 \nabla f_j $ is  the topological degree of the mapping
				
				\[\e\dfrac{\nabla f_j}{\|\nabla f_j\|}: S^{2n+1}_\e \to S^{2n+1}_\e ,\]
				for $ j = 1, \ldots,  p $.
				\item[(ii)] If $ n+1 $ is odd, then $ \chi(M_p) = 1 $. Moreover, one has $ \deg_0 \nabla f_i = 0 $ for $ i=1,2,\ldots, p. $
			\end{enumerate}
			
		\end{theorem}

		\vspace{0.2cm}

		\begin{theorem}\cite{RRE2}\label{tf}
			Let  $ F_{\eqref{eq: system prop}}=(f_1,f_2\dots, f_p) :(\mathbb{R}^{n+1},0) \to (\mathbb{R}^{p},0) $ be  a harmonic  first integral map for the system \eqref{eq: system prop} , 
			which has an isolated singularity at the origin. If any of the following conditions bellow occur, then the level sets of the solutions $u$ of  \eqref{eq: system prop} are included in the ball $ B^{n-p+1}\times B^{p-j} $ for any $ j=1,\dots, p $:
			\begin{enumerate}
				\item $  (n+1,p)=(4,2)  $ and $ \deg_0 \nabla f_1 = 0 $;
				\item $  (n+1,p)=(5,2) $ and $ \pi_1(M_p)=0 $, i.e., $ M_p $ is simply connected;
				\item $  (n+1,p)=(6,3)  $ and the manifold $  F^{-1}_{\eqref{eq: system prop}}(0) \cap S^n_\e $ is connected for small enough $ \e>0 $ or $ \deg_0 \nabla f_1 = 0 $;
				\item $  (n+1,p)=(8,5)  $ and the manifold $ F^{-1}_{\eqref{eq: system prop}}(0) \cap S^n_\e $ is not empty for small enough $ \e>0 $ or $ \deg_0 \nabla f_1 = 0 $;
				\item $ p=n $ and $ (n+1,p) \neq (4,3) $;
				\item $ n-p=1 $ and $(n+1,p) \neq (4,2) $;
				\item $  n-p+1=3   $ and $ (n+1,p)\neq (5,2) $, $  (n+1,p)\neq(6,3)  $  $  (n+1,p)\neq (8,5)  $. 
				
			\end{enumerate}
		\end{theorem}
		
		\vspace{0.2cm}

		The next two examples are based on the geometric/topological classifications presented in   \cite[section 4.2]{DDA} and can be naturally written for quasilinear system as \eqref{eq: system prop}.

		\vspace{0.2cm}
		
		\begin{example}
			Consider the germs of vector fields $X_1,X_2$ at $(\mathbb{R}^4,0)$ defined by
			\[
			X_1 = \frac{1}{2}(z^2-u^2)\frac{\partial}{\partial x} + zu \frac{\partial}{\partial y} -\frac{1}{2}(x_2-y^2)\frac{\partial}{\partial z} -xy\frac{\partial}{\partial u}
			\] 
			and 
			\[
			X_2 = zu\frac{\partial}{\partial x}+ \frac{1}{2}(z^2-u^2)\frac{\partial}{\partial y} -xy\frac{\partial}{\partial z} +\frac{1}{2}(x^2-y^2)\frac{\partial}{\partial u}.
			\]
			Since $X_1,X_2$ define an involutive distribution then it determine a  germ of codimension two foliation $\mathcal{F}$ at $(\mathbb{R}^4,0)$. Let $F = (f_1,f_2): (\bR^4,0) \to (\bR^2,0)$ be a map germ where $f_1=x^3-3 x y^2-3 z u^2+z^3$ and $f_2=-u^3+3uz^2+3x^2y-y^3$. Then, $df_j(X_i)\equiv 0, \forall i,j$.  Consequently, $F$ is a first integral map of $\mathcal{F}$. Note that $F$ has an isolated singularity at the origin. Furthermore, one has that $ \left\langle \nabla f_1, \nabla f_2 \right\rangle =0 $ and $ \|\nabla f_1 \|^2 =9u^4+18u^2z^2+9x^4+18x^2y^2+9y^4+9z^4= \|\nabla f_2 \|^2  $. Thus, $F$ is a horizontally weakly conformal map. Since $ F$ is polynomial, it follows from Theorem \ref{phm}  that $F$ is a harmonic morphism. Now, it follows from geometrical/topological description of the Milnor fibers  $M_{F}$, 
			that the leaves of $\mathcal{F}$ are all topologically equivalent to the torus minus three open discs removed.  In particular, it follows from  Theorem  \ref{tf} that  $ \deg_0 \nabla f_1 \neq 0 $.
		\end{example}
		
		\vspace{0.2cm}

		\begin{example}
			Let $\mathcal{F}$ be the codimension two germ of foliation at $(\mathbb{R}^4,0)$ determined by the vector fields 
			\[
			X_1 = \frac{5}{2}(z^4-6zû^2+u^4)\frac{\partial}{\partial x} + 10(z^3u-zu^3)\frac{\partial}{\partial y} -x\frac{\partial}{\partial z} -y\frac{\partial}{\partial u}
			\] 
			and 
			\[
			X_2 = 10(z^3u-zu^3)\frac{\partial}{\partial x}- \frac{5}{2}(z^4-6z^2u^2+u^4)\frac{\partial}{\partial y} -y\frac{\partial}{\partial z} +x\frac{\partial}{\partial u}.
			\]
			The foliation $\mathcal{F}$ admits a first integral map $F= (f_1,f_2): (\bR^4,0) \to (\bR^2,0)$
			with $f_1=5u^4z-10u^2z^3+z^5+x^2-y^2$ and $f_2=u^5-10u^3z^2+5uz^4+2xy$ which has an isolated singularity at the origin. Moreover, one has that $ \left\langle \nabla f_1, \nabla f_2 \right\rangle =0 $ and $ \|\nabla f_1 \|^2 =25u^8+100u^6z^2+150u^4z^4+100u^2z^6+25z^8+4x^2+4y^2= \|\nabla f_2 \|^2  $. Thus, $F$ is a horizontally weakly conformal map. Since $ F$ is polynomial, it follows from Theorem \ref{phm}  that it is a  harmonic morphism. Since the leaves of $\mathcal{F}$ are  topologically equivalent to a bitorus minus a open disc removed. Moreover,  
			one has that $ \deg_0 \nabla f_1 \neq 0 $. 
		\end{example}

		\subsubsection{The non-isolated singular case}

		Next, we extend \cite[Theorem 5]{RRE2} for the harmonic first integral maps case, with positive dimensional singular locus.

		\vspace{0.2cm}
		
		\begin{proposition} Let  $ F_{\eqref{eq: system prop}}=(f_1,f_2\dots, f_p) :(\mathbb{R}^{n+1},0) \to (\mathbb{R}^{p},0) $ be  a harmonic  first integral map with $\dim\,\Sing F_{\eqref{eq: system prop}} \ge 0 $. 
			The manifold $  M_j $ is homeomorphic to $  M_p\times {B}^{p-j}  $, for any $ j=1,\dots,p $.
		\end{proposition}
		
		\begin{proof}
			Since  harmonic morphisms satisfy the condition \eqref{eq:main} and has an isolated critical value at the origin, the result follows directly from \cite[Theorem 6.3 ]{DA} applied successively .
		\end{proof}
		
		\vspace{0.2cm}
		
		The next result shows that the Euler characteristic of each  integral manifold is equal to the Euler characteristic of the its intersections.
		
		\vspace{0.2cm}

		\begin{corollary}
			Let  $ F_{\eqref{eq: system prop}}=(f_1,f_2\dots, f_p) :(\mathbb{R}^{n+1},0) \to (\mathbb{R}^{p},0) $ be  a harmonic first integral map, with $p \geq 2,$  and $\dim\,\Sing F_{\eqref{eq: system prop}} \ge 0 $. Then $ \chi(M_p) =\chi(M_j) $, for any $ j=1,\dots,p $. In particular, $ \chi(N_j) = \chi(M_p) $, for any $ j=1,\dots,p $.
		\end{corollary}
		\begin{proof}
			It follows from \cite[Corollary 6.4, Corollary 6.5 and Corollary 6.6]{DA}.
		\end{proof}

		\vspace{0.2cm}
		
		\begin{remark}
			Let us point out that for any $ j=1,\dots,p $, the first integral $f_j :(\mathbb{R}^{n+1},0) \to (\mathbb{R},0) $ of the system $ \eqref{eq: system prop} $ admits two Milnor fibers, namely: $ N_j=f^{-1}(\delta)\cap B^{n+1}_\epsilon $ and  $ N_j^{-}:=f^{-1}(-\delta)\cap B^{n+1}_\epsilon $, where $ 0<\delta \ll \epsilon $. Therefore, one has $ \chi(N_j^{-})=\chi(N_j)=\chi(M_p) $. For more details, see \cite[Corollary 6.4]{DA}
		\end{remark} 
		
	

	\vspace{0.2cm}
	
	As we have seen, the existence of bundle associated with  singularities of a harmonic first integral map,  make it possible to study the topology of manifolds without singularities near zero $M_j$,   which are ``close'' of the analytic sets $V_{F^j_{\eqref{eq: system prop}}}:=f_1^{-1}(0)\cap f_2^{-1}(0)\cap \dots \cap f_j^{-1}(0)$,  with $ 1\le j \le p-1 $, and $ V_{F_{\eqref{eq: system prop}}}:=F_{\eqref{eq: system prop}}^{-1}(0)$ for $ j=p $. In this sense, we apply one more result of \cite{DA} which provides a criterion for computing the Euler characteristic of $M_j$.
	
	\vspace{0.2cm}
	
	\begin{proposition}\cite[Theorem 7.4]{DA}
		Let  $ F_{\eqref{eq: system prop}}=(f_1,f_2\dots, f_p) :(\mathbb{R}^{n+1},0) \to (\mathbb{R}^{p},0) $ be  a harmonic first integral map, with $p \geq 2,$  and $\dim\,\Sing F_{\eqref{eq: system prop}} \ge 0 $. Then:
		
		\begin{enumerate} 
			\item $\chi(V_{F^j_{\eqref{eq: system prop}}} \cap S^n_\e ) = 2\chi(M_p) $, if $ n+1 $ is even and $ j \in \{1,\dots, p\} $ is odd, for any small enough $ \e>0 $.
			
			\item $\chi(V_{F^j_{\eqref{eq: system prop}}} \cap S^n_\e ) = 2-2\chi(M_p) $, if $ n+1 $ is odd and $ j \in \{1,\dots, p\} $ is odd, for any small enough $ \e>0 $.
		\end{enumerate}
		
	\end{proposition}

	\section{Geometry of harmonic first integral maps}
	
	The authors of \cite{RADG}  generalized  the main result of \cite{PY}, which ensure, about special hypotheses over a germ of  harmonic morphism, that the Milnor sphere fibration is itself a harmonic morphism.  In  this section, we will make use of this important fact  to describe  the  geometric properties of Milnor fibers of the harmonic first integral maps. 
	
	\vspace{0.2cm}
	
	Our main results are as follows.
	
	\vspace{0.2cm}
	
	\begin{proposition}\label{p1}
		Let  $ F_{\eqref{eq: system prop}}=(f_1,f_2\dots, f_p) :(\mathbb{R}^{n+1},0) \to (\mathbb{R}^{p},0) $, be  a harmonic  first integrals map such that $ f_1,f_2\dots, f_p $ are homogeneous polynomials of the same degree $ d $ and $ p\ge 3 $.  If $  V_{F_{\eqref{eq: system prop}}} = \{0\} $,  then the level sets of the solutions $u$ of  \eqref{eq: system prop} are included in a Hopf fibrations fiber. In particular,  one has that $n= 3,7,15  $ ,  $ p= 3, 5, 9 $ and $ d=2 $.
	\end{proposition}
	
	\vspace{0.2cm}

	\begin{theorem}\label{t2}
		Let  $ F_{\eqref{eq: system prop}}=(f_1,f_2\dots, f_p) :(\mathbb{R}^{2n},0) \to (\mathbb{R}^{2},0) $, be  a harmonic  first integral map such that $ f_1,f_2\dots, f_p $ are homogeneous polynomials of the same degree $ 2 $. If  $ \dim \,V_{F_{\eqref{eq: system prop}}} \ge 0, $  then $ \Psi_{F_{\eqref{eq: system prop}}} $ is equivalent for $ \Psi_f $, where $ f:(\bC^n,0) \to (\bC) $ is given by $ f(z_1, \dots, z_n)=\lambda_1 z_1^2 + \dots +\lambda_n z_n^2$ with $ \lambda_k \ge 0 $ for any $ k=1,\dots,n $ and not all zero. In particular, one has that $ \Sing F_{\eqref{eq: system prop}} =\{0\} $, and the level sets of the solutions $u$ of  \eqref{eq: system prop} are included in a  $ (n-1) $-sphere.
	\end{theorem}
	
	\vspace{0.2cm}
	
	The proof of Theorem \ref{t2} and Proposition \ref{p1} follows from the next result, which shows that under special conditions, Milnor fibrations associated with,  in a certain sense, equivalent harmonic morphisms, have the same topological and geometric structure.
	
	\vspace{0.2cm}

	\begin{lemma}\label{l1} Let  $ F, G :(\mathbb{R}^{n+1},0) \to (\mathbb{R}^{p},0) $ be  harmonic morphisms. If there is an isometry $ \phi \in O(\bR^{n+1}) $ such that $ F=  G\circ \phi  $ then the Milnor fibrations $ \Psi_F $ and $ \Psi_G $ are equivalents. 
		
	\end{lemma}
	\begin{proof}
		It follows from hypothesis that $ \phi(V_F)=V_G$. Consequently, $V_F$ and $ V_G $ are diffeomorphic, in particular, since $ \phi $ is a isometry, for any small enough $ \e>0 $, one has $ \phi\left( S^n_\e \m \left( V_F \cap S^n_\e\right) \right)  = S^n_\e \m \left( V_G \cap S^n_\e\right)  $. Consider the Milnor fibrations $ \Psi_G $ and $ \Psi_F $, one has that   $ {\Psi_G}(x)=\Psi_F \circ \phi(x) $ for any $ x\in S^n_\e \m \left( V_F \cap S^n_\e\right) $. Therefore, one has that $ \Psi_F = \Psi_G \circ \phi_{|_{S^n_\e \m \left( V_F \cap S^n_\e\right)}} $, and the fibrations $ \Psi_F $ and $ \Psi_G $ are equivalents. 
	\end{proof}

	\vspace{0.2cm}

	\begin{proof}[\textbf{Proposition \ref{p1}}]
		Since $ f_1,f_2\dots, f_p $ are homogeneous polynomials of the same degree $ d $ and  $  V_{F_{\eqref{eq: system prop}}} = \{0\} $, it follows from 
		\cite[Lemma 1]{PY}  that $\Psi_{F_{\eqref{eq: system prop}}}: S^n_\e \to S^{p-1}_\eta$ is a full submersive harmonic morphism wich is non-constant map given by the restriction of $ F_{\eqref{eq: system prop}}$.  Consequently,  it follows from \cite[Theorem 5.6.5]{BW} that $\Psi_{F_{\eqref{eq: system prop}}}$   there exist an isometry $ \phi \in O(\bR^{n+1}) $ such that $ F_{\eqref{eq: system prop}}=  f\circ \phi  $, where $ f $ denotes one of  Hopf fibrations $ S^3_\e  \to S^2_\eta $, $ S^7_\e  \to S^4 _\eta$, $ S^{15}_\e  \to S^8 _\eta$. Now the result follows from the Lemma \ref{l1}.
	\end{proof}
	
	\vspace{0.2cm}

	\begin{proof}[\textbf{Theorem \ref{t2}}]
		Since  $ f_1,f_2\dots, f_p $ are homogeneous polynomials of the same degree $ 2 $,   It follows from \cite[Theorem 5.11]{YO} that exist $ \lambda_1, \dots, \lambda_n \ge 0 $  not all zero, and an isometry $ \phi \in O(\bR^{2n}) $ such that $ F_{\eqref{eq: system prop}}=  f\circ \phi  $, where $ f:(\bC^n,0) \to (\bC) $ is given by $ f(z_1, \dots, z_n)=\lambda_1 z_1^2 + \dots +\lambda_n z_n^2$. One has that $ \Sing f = \{0\} $ and $ \phi $ is a diffeomorphism,  thus $ \Sing F = \{0\} $. Next, it is follows from \cite[Theorem 6.5 and Theorem 7.2 ]{Mi} that  the fiber $ \Psi_f^{-1}(\theta) $, with $ \theta \in S_\eta^1 $,  has the same homotopy type of a bouquet of $(n-1)$-dimensional spheres $ \bigvee_{i = 1}^{\mu(f)} S^{n-1}_{i}$, where $\mu(f) $ is the Milnor number of the singularity. Since $\mu(f) =1 $,  the Lemma \ref{l1} ensure that the level sets of the solutions $u$ of  \eqref{eq: system prop} are included in $ S^{n-1}$, up to isometries.
		
	\end{proof}
	
	\vspace{0.2cm}

	\begin{theorem}\label{tf2} Let  $ F_{\eqref{eq: system prop}}=(f_1,f_2\dots, f_p) :(\mathbb{R}^{n+1},0) \to (\mathbb{R}^{p},0) $, be  a harmonic  first integrals map such that $ f_1,f_2\dots, f_p $ are homogeneous polynomials of the same degree $ d $. If  $ \dim \,V_{F_{\eqref{eq: system prop}}} \ge 0, $ and one of the conditions below occurs:
		
		\begin{enumerate}
			\item $ p=3 $
			
			\vspace{0.2cm}
			
			\item $ p\ge4 $ and  $\Psi_{F_{\eqref{eq: system prop}}}$ is horizontally homothetic\footnote{See \cite[Definition 2.4.18]{BW}.}.
		\end{enumerate}
		Then the fibers of the Milnor sphere fibration $\Psi_{F_{\eqref{eq: system prop}}}$, are minimal.
	\end{theorem}
	
	\begin{proof}
		It is follows from  \cite[Corollary 5.2]{RADG}  that $\Psi_{F_{\eqref{eq: system prop}}}$ is a full submersive harmonic morphism.	Next, if $ p=3 $ the result follows from \cite[Theorem 4.5.4]{BW}. If $ p\ge 4 $ and  $\Psi_{F_{\eqref{eq: system prop}}}$ is horizontally homothetic, the result follows from \cite[Theorem 4.10]{BW}.
		
	\end{proof}
	
	\vspace{0.2cm}

	\begin{example} 	Let $\mathcal{F}$ be the dimension five germ of foliation at $(\mathbb{R}^8,0)$ determined by $$\omega =\omega_1 \wedge \omega_2 \wedge \omega_3,$$ where each  $\omega_j$ is a 1-form determined by
		\[
		\omega_1 = 2x_1 dx_1 + 2x_2dx_2 - 2x_3dx_3 -2x_4dx_4 +3x_5dx_5+3x_6dx_6-3x_7dx_7-3x_8dx_8, 
		\]
		\[
		\omega_2 = 2x_3dx_1-2x_4dx_2+2x_1dx_3-2x_2dx_4+3x_7dx_5-3x_8dx_6+3x_5dx_7-3x_6dx_8,
		\]
		and 
		\[
		\omega_3 = -2x_4dx_1+2x_3dx_2+2x_2dx_3-2x_1dx_4+3x_8dx_5-3x_7dx_6-3x_6dx_7+3x_5dx_8.
		\]
		Then, $\mathcal{F}$ admits a first integrals map $ F=(f_1,f_2,f_3):(\mathbb{R}^8,0) \to (\mathbb{R}^3,0) $, where
		
		$$ \left\lbrace \begin{array}{cl}
			f_1(x,y,z,w,a,b,c,d)&=2(x^2+y^2)+3(a^2+b^2) -2(z^2+w^2)-3(c^2+d^2)     \medskip\\
			f_2(x,y,z,w,a,b,c,d)&= 6ac-6bd-4wy+4xz  \medskip\\
			f_3(x,y,z,w,a,b,c,d)&= 6ad+6bc+4wx+4yz  \medskip\\

		\end{array}\right. $$ 
		
		By a direct calculation, one has $ \left\langle \nabla f_1, \nabla f_2 \right\rangle = \left\langle \nabla f_1, \nabla f_3 \right\rangle = \left\langle \nabla f_2, \nabla f_3 \right\rangle =0 $ and $ \|\nabla f_1 \|^2 	= \|\nabla f_2 \|^2=\|\nabla f_3 \|^2 = 36a^2+36b^2+36c^2+36d^2+16w^2+16x^2+16y^2+16z^2$. Thus, $ F$ is a horizontally weakly conformal map.  Since $ F$ is polynomial, it follows from Theorem \ref{phm}  that $ F$ is a  harmonic first integral map. Now,  from Theorem \ref{tf2}, the fibers of the  Milnor sphere fibration $ \Psi_{F_{\eqref{eq: system prop}}}: S^7 \m K_\e \to S^2 $ are minimal.

	\end{example}

	\vspace{0.2cm}

	\begin{example}\cite{PY} Let $\mathcal{F}$ be the dimension four germ of foliation at $(\mathbb{R}^8,0)$ determined by 
		\[
		\omega =\omega_1 \wedge \omega_2 \wedge \omega_3 \wedge \omega_4
		\]
		where each $\omega_j$ is a $1$-form given by
		\[
		\omega_1 = x_5dx_1 - x_6dx_2 -x_7dx_3 -x_8dx_4+x_1dx_5 -x_2dx_6 -x_3dx_7 -x_4dx_8,
		\] 	
		\[
		\omega_2 = -x_6dx_1 + x_5dx_2 - x_8dx_3 -x_7dx_4 + x_2dx_5 -x_1dx_6 -x_4dx_7 -x_3dx_8,
		\]	
		\[
		\omega_3 = x_7dx_1 -x_8dx_2 + x_5dx_3 +x_6dx_4 +x_3 dx_5 +x_4dx_6 +x_1dx_7 -x_2dx_8
		\]	
		and
		\[
		\omega_4 = -x_8dx_1-x_7dx_2 -x_6dx_3 +x_5dx_4+x_4dx_5 - x_3dx_6 -x_2dx_7 -x_1dx_8.
		\]	
		Hence, $\mathcal{F}$ admits a first integrals map $ F=(f_1,f_2,f_3, f_4):(\mathbb{R}^8,0) \to (\mathbb{R}^4,0) $, where
		
		$$ \left\lbrace \begin{array}{cl}
			f_1(x,y,z,w,a,b,c,d)&=ax-by-cz-dw    \medskip\\
			f_2(x,y,z,w,a,b,c,d)&= ay+bx-cw+dz \medskip\\
			f_3(x,y,z,w,a,b,c,d)&= az+bw+cx-dy \medskip\\
			f_4(x,y,z,w,a,b,c,d)&= aw-bz+cy+dx \medskip\\

		\end{array}\right. $$ 
		
		By a straightforward computation one has $ \left\langle \nabla f_i, \nabla f_j \right\rangle = 0 $ for $ i,j=1,2,3,4 $ and $ i\neq j $. Moreover,  $ \|\nabla f_i \|^2 = a^2+b^2+c^2+d^2+w^2+x^2+y^2+z^2$, for any $ i=1,2,3,4 $. Thus, $ F$ is a a horizontally weakly conformal map.  Since $ F$ is polynomial, it follows from Theorem \ref{phm}  that it is a  harmonic first integral map.   Moreover, $ \Psi_{\phi} $ is a horizontally homothetic. Therefore, it follows from Theorem \ref{tf2}  that the fibers of the  Milnor sphere fibration $ \Psi_{F_{\eqref{eq: system prop}}}: S^7 \m K_\e \to S^2 $ are minimal.

	\end{example}

\end{document}